\documentstyle[12pt]{article}
\renewcommand{\baselinestretch}{1.3}

\newtheorem {th}{Theorem}[section]
\newtheorem {lem}[th]{Lemma}

\def\Cox{\hfill \Box}
\def\sym{\bigtriangleup}
\def\dj{\mbox{ is disjoint from }}

\def\dconv{\, {\stackrel {{\cal D}} {\rightarrow}}}

\def\sf{\sigma\mbox{-field}}

\def\ee{\epsilon}
\def\E{{\bf{E}}}
\def\P{{\bf{P}}}
\def\N{\hbox{I\kern-.2em\hbox{N}}}
\def\R{\hbox{I\kern-.2em\hbox{R}}}
\def\Z{{\bf{Z}}}

\def\C{{\cal{C}}}

\def\F{{\cal{F}}}
\def\|{\, | \, }
\def\one{{\bf 1}}
\def\Tree{{\bf T}}

\begin{document}

\begin{titlepage}
\begin{center}
{\large \bf CHOOSING A SPANNING TREE FOR THE INTEGER LATTICE
UNIFORMLY}  \\ 
Running Head: RANDOM SPANNING TREES
\end{center}
\vspace{5ex}
\begin{flushright}
Robin Pemantle \footnote{This research supported by a National 
Science Foundation postdoctoral fellowship and by a Mathematical
Sciences Institute postdoctoral fellowship.  }\\
Dept. of Mathematics, White Hall \footnote{Now at the Department
of Mathematics, Oregon State University, Corvallis, OR 97331}\\
Cornell University \\
Ithaca, NY 14853 \\  
\today 

\end{flushright}

\vfill

{\bf ABSTRACT:} \break
Consider the nearest neighbor graph for the integer lattice $\Z^d$ in 
$d$ dimensions.  For a large finite piece of it, consider choosing a
spanning tree for that piece uniformly among all possible subgraphs
that are spanning trees.  As the piece gets larger, this approaches
a limiting measure on the set of spanning graphs for $\Z^d$.  This
is shown to be a tree if and only if $d \leq 4$.  In this case, 
the tree has only one topological end, i.e. there are no doubly 
infinite paths.  When $d \geq 5$ the spanning forest has infinitely
many components almost surely, with each component having one or
two topological ends.

\vfill



\end{titlepage}

\section{Introduction} \label{intro}

Let $\Z^d$ be the nearest neighbor graph on the $d$-dimensional 
integer lattice, so there is an edge between $(v_1, \ldots v_d)$
and $(w_1, \ldots , w_d)$ if and only if $\sum_i |v_i - w_i | = 1$.
The term {\em subgraph} will be used to denote any subcollection
of these edges.  A subgraph of $\Z^d$ {\em spans} $\Z^d$ if it
contains at least one edge incident to each vertex.  A graph
is a {\em forest} if it has no loops and a {\em tree} if it
is a connected forest.  A spanning tree on $\Z^d$ is thus a
connected, loopless subgraph of $\Z^d$ that spans $\Z^d$.

For measure theoretic purposes, subgraphs are viewed as
maps from the set of edges of $\Z^d$ to $\{ 0,1 \}$.
Topologize the space of all subgraphs by the product topology,
generated by the {\em cylinder sets}, namely those sets 
depending on only finitely many edges.  There is a Borel $\sf$
for this topology and it is also generated by the 
{\em elementary} cylinder sets, $C(A)$, where $A$ is a
finite set of edges and $C(A)$ is the set of subgraphs 
containing all the edges in $A$.  For measures on the Borel $\sf$,
$\nu_n \rightarrow \nu$ weakly iff $\nu_n (C) \rightarrow \nu (C)$
for every cylinder set $C$; it suffices to check this for elementary
cylinder events $C(A)$. 

This paper is concerned with the following method of
picking a spanning tree on $\Z^d$ at random.  Let $B_n$ be
the box of diameter $2n$ centered at the origin, so it has
$(2n+1)^d$ vertices and all the nearest neighbor edges between
these vertices.  Let $|v-w|$ denote the metric $\max \{ v_i - w_i \}$;
this is convenient for counting and for making $B_n$ a sphere,
although any equivalent metric could be substituted throughout
with no change to the theorems.  There are finitely many spanning 
trees on $B_n$ so there is a uniform measure $\mu_1 (B_n ) $ on 
spanning trees of $B_n$.  Any spanning tree on $B_n$ is a
subgraph of $\Z^d$ so one may view the measure $\mu_1 (B_n )$
as a measure on subgraphs of $\Z^d$.  It turns out that these
measures converge weakly as $n \rightarrow \infty$ to a measure
$\mu$ on spanning forests of $\Z^d$.  For notational
convenience, abbreviate $\mu \{ \Tree \, :\, \cdots \}$ to 
$\mu ( \cdots )$.  

The main tool for proving this basic result is the
equivalence (for finite graphs) between uniform spanning
trees and random walks.  Together with the further equivalence
between random walks and electrical networks, this provides
a basis for proving that the measures $\mu_1 (B_n)$ converge as
well as proving some ergodic properties of the limiting measure
$\mu$ that will be important later.  This groundwork is 
laid in section~\ref{equivalences}.  

The rest of the paper is concerned with the geometry
of the typical sample from the measure $\mu$.  It is easy to
see that $\mu$ concentrates on spanning forests of $\Z^d$.  The
first result is that in dimensions $d \leq 4$ the measure 
concentrates on spanning trees, while in dimensions $d \geq 5$,
the spanning forest will almost surely have infinitely many 
components.  The shape can be further described by the number
of topological ends.  For a tree, the number of topological ends
is just the number of infinite, self-avoiding paths from any 
fixed vertex.  It turns out that when $d \leq 4$ the measure
concentrates on spanning trees with only one end.  
When $d \geq 5$ the measure concentrates on spanning forests
in which each of whose components has one or two topological ends.

The machinery used to prove these shape results is 
Lawler's theory of loop-erased random walks (LERW).  These are 
defined in section~\ref{LERW} and the required basic results 
about LERW are referenced or proved.  The shape results
are then proved in section~\ref{shape}.  

\noindent{{\bf Acknowledgement:}} All of the questions studied in 
this paper were asked by Russ Lyons.

\section{Uniform spanning trees, random walks and electrical networks}
\label{equivalences}

For any connected finite graph $G$, let $\mu_1 (G)$ be the 
uniform measure on spanning trees of $G$, as in section~\ref{intro}.
Let $v$ be any vertex of $G$.  The following defines a measure 
$\mu_2 (G,v)$ which will turn out to be the same as $\mu_1 (G)$,
independently of $v$.  Let $\gamma = \gamma (0), \gamma (1), 
\ldots$ be a simple
random walk (SRW) on $G$ starting from $v = \gamma (0)$.  Let 
$\Tree (\gamma )$ be the subgraph of $G$ containing precisely
those edges $\overline {\gamma (i) \, \gamma (i+1)}$ for which there is
no $j < i$ with $\gamma (j) = \gamma (i+1)$.  Another way to describe
$\Tree (\gamma )$ is ``walk along gamma and draw in each edge
as you go except when drawing in an edge would close a loop''.
The graph $\Tree (\gamma )$ depends only on $\gamma (0) , \ldots 
, \gamma (\tau )$ where $\tau$ is the first time $\gamma$ has
visited every vertex.
The SRW measure on paths $\gamma$ projects to a measure
$\mu_2 (G,v)$ on subgraphs of $G$.  By viewing these edges as
oriented from $\gamma (i)$ to $\gamma (i+1)$ is is easy to see that
the resulting subgraph is a spanning tree on $G$ oriented away
from $v$.
 
\begin{lem} \label{unif=srw}
For any vertex $v$ of a finite graph $G$, $\mu_1 (G) = \mu_2 (G,v)$.
\end{lem}
 
\noindent{Proof:}  This result is due to Diaconis and Doyle;
a more complete account can be found in Aldous (1988) or 
Broder (1988).  Let $\{ v_i \,:\, i \in \Z \}$ be the
stationary Markov chain corresponding to SRW on $G$.  Let
$\Tree_i$ be the rooted tree whose oriented edges are just
those edges $\overline {v_j \, v_{j+1}}$ for which $v_{j+1}$
is distinct from every $v_k$ for $i \leq k < j$.  It is easy to
check that $\Tree_i$ is indeed loopless and almost surely
connected and that all edges are oriented away from $v_i$, 
which is taken to be the root.  Furthermore, it may be
verified that $\{ \Tree_i \}$ is a stationary Markov chain
on the space of rooted spanning trees of $G$ and that 
a unique stationary measure for it is given
by letting the measure of each rooted tree be proportional
to the number of neighbors of the root.  This means that 
conditioning on the root of the tree (which is just $v_0$)
leaves a uniform unrooted spanning tree.  Now the SRW measure from 
$v$ is just the stationary Markov measure on $\{ v_i \,:\,
i \geq 0 \}$ conditioned on $v_0 = v$.  Thus $\mu_2 (G,v)$ is
distributed as $\Tree (v_0 , v_1 , \ldots )$, where $\{ v_i \}$ 
are a stationary Markov chain conditioned on $v_0 = v$.
This has just been shown to be uniform, and the proof is done.
$\Cox$

For any edge $e = \overline {v \, w}$ in a finite 
graph $G$, define the {\em contraction} of $G$ by $e$ to be the
graph $G / e$ gotten by removing $e$ and identifying $v$ and $w$.
This may result in parallel edges, which must still be
regarded as distinct, or in loops (edges whose endpoints 
are not distinct) which may for the purposes
of what follows be thrown away.  The {\em deletion} of $e$
from $G$ is just the graph $G-e$ consisting of all edges of $G$ except 
$e$.  Contraction commutes and associates with deletion, 
so it makes sense to speak of the graph $G$ with $e_1,
\ldots e_r$ contracted and $e_1', \ldots , e_s'$ deleted.
Note that there are natural identifications $\phi^{(-e)}$
and $\phi^{(/e)}$ between edges of
$G$ other than $e$ and edges of either $G/e$ or $G-e$.

Now another measure will be defined on subgraphs of a given graph
$G$ that turns out to be the same as $\mu_1 (G)$.  
Let $\C = e_1 , e_2, \ldots$ be any enumeration of the edges
of a finite graph $G$.  Define $\mu_3 (G , \C)$
recursively as follows.  The start of the recursion is that
if $G$ is a single vertex then $\mu_3 (G)$ is the pointmass at $G$.
To continue the recursion, assume that $\mu_3 (G)$ is defined
for all contractions and deletions of $G$ and all enumerations.
To define $\mu_3 (G,\C )$ begin by throwing out all loops and
putting a 1 ohm resistor along each edge.  Put the terminals
of a battery at the two ends of $e_1$.  Look at the total
current that flows through the battery and see what fraction
of it flows through the resistor at $e_1$.  Call this fraction
$p$.  There is a random walk interpretation for $p$: it is the
probability that a simple random walk started at one end of the edge 
$e$ reaches the other end for the first time by moving along $e$.
Let the $\mu_3 (G)$ measure give probability $p$ to
the event $e_1 \in \Tree$ and $1-p$ to the complementary
event.  The specification of $\mu_3$ is completed by stating
the conditional distributions of $\mu_3$ given $e_1 \notin \Tree$
and $e_1 \in \Tree$.  To do this write $\C' = e_2 , e_3 , \ldots$,
where $e_2 , e_3 , \ldots$ are viewed as edges in $G - e$ or
$G/e$ via the natural identifications $\phi^{(-e)}$ and
$\phi^{(/e)}$.  Then the distribution of $\mu_3 (G,\C)$
given $e_1 \notin \Tree$ is just $\mu_3 (G-e_1 , \C')$, which is
a measure on subgraphs of $G-e_1$, hence on subgraphs of $G$
via $\phi^{(-e)}$.
Let the distribution of $\mu_3 (G,\C )$ given $e_1 \in \Tree$
be given by adding the edge $e_1$ to a subgraph of $G$
chosen by picking a subgraph of $G/e_1$ from $\mu (G/e_1 , \C')$
and viewing it as a subgraph of $G$ by the natural identification
$\phi^{(/e)}$.

\begin{lem} \label{electric}
For any enumeration $\C = e_1 , e_2 , \ldots$ of the edges of
a finite connected graph $G$, the measure $\mu_3 (G,\C )$
is equal to $\mu_1 (G)$.
\end{lem}

\noindent{Proof:}  The idea of the proof is that $\mu_1$ 
satisfies the same recursion as $\mu_3$.  Begin by observing
that the spanning trees of $G$ that do not contain an edge
$e$ are in one to one correspondence with the spanning trees
of $G-e$.  Secondly, observe that the
spanning trees of $G$ that do contain $e$ are in one to one
correspondence with the spanning trees of $G/e$, where the
correspondence is given by subtracting the edge $e$.  This
is because the identification of the endpoints of $e$ in
$G/e$ makes a set of edges of $G/e$ contain a loop
if and only if the set together with $e$ contains a loop in $G$.
It is clear that single edge loops of $G/e$ may be thrown out.

These observations imply that $\mu_1 (G)$
conditioned on $e \in \Tree$ is just $\mu_1 (G/e)$ and
$\mu_1 (G)$ conditioned on $e \notin \Tree$ is just 
$\mu_1 (G-e)$.  The next thing is to see that the event
$B = \{ e_1 \in \Tree \}$ has the same probability under $\mu_1$
as it does under $\mu_3 (G,\C )$ for any enumeration $\C$ 
beginning with $e_1 = \overline {v \, w}$.  
By Lemma~\ref{unif=srw}, $\mu_1 (B)$
is the probability that a SRW on $G$ from $v$ has just traveled
across $e$ when it hits $w$ for the first time.  By the well-known
correspondence between random walks and electrical networks
(see Doyle and Snell, section 3.4), this is precisely the
fraction $p$ of the current that flows across $e_1$ in the
electrical scenario used to define $\mu_3$.  

Now it follows that if $\mu_1 (G/e_1) = \mu_3 (G/e_1 , \C' )$ 
and if either $G-e$ is disconnected or 
$\mu_1 (G-e_1) = \mu_3 (G-e_1 , \C' )$, then $\mu_1 (G) 
= \mu_3 (G, \C )$.  The initial conditions are certainly
the same: if $G$ is a single vertex then $\mu_1 (G)$ is the pointmass
at $G$.  By induction on the number of edges, it follows that
$\mu_1 (G) = \mu_3 (G, \C )$ for all finite connected graphs 
and enumerations.    $\Cox$

\begin{th} \label{convergence}
Let $\{ B_n \}$ be any sequence of finite sets of edges of $\Z^d$,
$d \geq 2$, converging to $\Z^d$ in the sense that any edge is 
in all but finitely many sets $B_n$.  Then the measures 
$\mu_1 (B_n)$ converge weakly to a limiting measure $\mu$
in the sense that $\mu_1 (B_n ) (C) \rightarrow \mu (C)$ for
any cylinder event $C$.  The measure $\mu$ is concentrated on
spanning forests of $\Z^d$ and is translation invariant.
\end{th}

\noindent{Proof:}  For weak convergence it suffices to show
that $\mu_1 (B_n ) (C)$ converges for the special
case where $C$ is the event $C(A)$ that all edges in a 
finite set $A$ are in the random subgraph.  This is because
the probabilities of $C(A)$ determine the probabilities of
all cylinder events by inclusion-exclusion, and because if
all cylinder probabilities converge the limits of these must
define a measure.  

Proceed by fixing a set $A = e_1 , \ldots , e_k$.  When
$n$ is sufficiently large so $A \subseteq B_n$, 
let $\C_n$ be an enumeration of the edges of $B_n$ that begins 
with $e_1, \ldots , e_k$.  Then by the previous Lemma,
$\mu_1 (B_n )(C(A)) = \mu_3 (B_n , \C ) (C(A)) = 
\prod_{j=1}^k p_j^{(n)}$ where $p_j^{(n)}$ is the 
$\mu_3 (B_n /e_1 / \cdots
/ e_{j-1} , \C^{(j-1)})$ probability of $\{ e_j \in \Tree \}$.
This is just the fraction of current that flows through
$e_j$ when a battery is placed across $e_j$ in the resistor
network $B_n /e_1 / \cdots / e_{j-1}$.  

Consider for a moment the special case where $B_n$ is a box of
diameter $2n$ centered at the origin.  Then for $r>0$, $B_n$ is
just $B_{n+r}$ with a lot of edges removed.  Since
contraction and deletion commute, $B_n /e_1 / \cdots / e_{j-1}$
is just a deletion of $B_{n+r} /e_1 / \cdots / e_{j-1}$.
It follows from Raleigh's Monotonicity Law (Doyle and 
Snell Chapter 4) that more current flows in $B_{n+r} /e_1 / \cdots 
/ e_{j-1}$ than in $B_n /e_1 / \cdots / e_{j-1}$.  Since the 
same current flows directly across the edge $e_j$, it follows
that $p_j^{(n)} \geq p_j^{(n+r)}$ and by taking the product
that $\mu_1 (B_n ) (C(A)) \geq \mu_1 (B_{n+r}) (C(A))$.
The sequence of probabilities is therefore decreasing in 
$n$ and must converge for each $A$.  

For general $B_n$, note
that the $B_n$ eventually contain any finite box and are
each contained in some finite box.  The monotonicity proof
worked for any graphs, one containing the other.  Then the
probabilities $\mu_1 (B_n) (C(A))$ interlace the sequence
of probabilities of $C(A)$ for boxes of diameter $2n$ and
hence converge to the same limit.  

The rest is immediate.  There are no loops in the final
measure $\mu$, because any loop $e_1 , e_2 , \ldots , e_k$ is a finite
cylinder event and has probability zero under each $\mu_1 (B_n)$.
Also, the event that vertices $v_1 , \ldots , v_k$ are a component
not connected to the rest of the graph is a cylinder event on
any box $B_n$ big enough to contain all edges incident to any $v_i$.
The $\mu_1 (B_n)$ probability of this event is zero, since $\mu_1 (B_n)$
concentrates on connected graphs, so the limit is zero.  
For stationarity, note that $\mu_1 (B_n ) (C(\pi A))
= \mu_1 (\pi^{-1} B_n ) (C(A))$ for any translation $\pi$.
The interlacing argument shows that using the sequence
$\pi^{-1} B_n$ in place of $B_n$ does not affect the limit,
so $\mu (C(\pi A)) = \mu (A)$ for any event $C(A)$.  These events
determine the measure, hence $\mu$ is translation invariant.
$\Cox$

For any set $A$ of edges, let $\sigma (A)$ denote as usual
the $\sf$ generated by the events $C(A')$ for $A'$ a 
finite subset of $A$.  Let $\F$ denote the tail $\sf$, 
which is just the intersection of $\sigma (A)$ over all cofinite
sets $A$.

\begin{th} \label{tail}
Let $\mu$ be the measure defined above on spanning forests
of $\Z^d$.  Then the tail field is trivial, i.e. $\mu (C)
= 0$ or $1$ for every $c \in \F$.
\end{th}

\noindent{Proof:}  First the electrical viewpoint will be used
to reduce the statement to a more specialized proposition
and then the random walk construction will be used to prove
the proposition.  

Begin with the device used to prove Kolmogorov's zero-one
law: an event is trivial if it is independent from
every event in a sufficiently large set.  Letting $C$ be any tail 
event, it suffices to show that if $\mu (C) > 0$ then the conditional
probabilities $\mu ( . \| C)$ agree with $\mu$ on elementary
cylinder sets.  For $n > 0$, let $B_n$ be
boxes of diameter $2n$ centered at the origin and let $C_n$ be
cylinder sets in $\sigma (\Z^d \setminus B_n )$ such that
$\mu (C_n \sym C) \rightarrow 0$.  In particular, the sequence
$\{ \mu (C_n) \}$ has a positive $\lim\inf$ and it will suffice 
to show that for
each finite set of edges $A$, $\mu (C(A) \| C_n ) 
\rightarrow \mu (C(A))$ for $n$ such that $\mu (C_n) \neq 0$.  
By Lemma~\ref{electric}, it suffices
to show that for any sequence of boxes $B_n'$ big enough so
that $C_n \in \sigma (B_n')$, $\,\mu_3 (B_n') (C(A)) \|
C_n) \rightarrow \mu_ 3 (B_n') (A)$ at least for those
$n$ such that $\mu_3 (B_n') (C_n ) \neq 0$.  

To do this, consider the electrical networks $G_1$
and $G_2$ where $G_1$ is just $B_n$ and $G_2$ is gotten by 
contracting all edges outside of $B_n$, which is electrically
the same as short circuiting the boundary, $\partial B_n$, of the 
box $B_n$.  I claim that $\mu_3 (B_n') (C(A) \| D)$ is bounded below by
$\mu_3 (G_2 ) (C(A))$ and above by $\mu_3 (G_1 ) (C(A))$
for any event $D \in \sigma (B_n' \setminus B_n )$.
To see this, let $\C$ be an enumeration of the edges in
$B_n'$ beginning with those not in $B_n$.  The event $D$
is a union of cylinder events that specify precisely
which edges in $B_n' \setminus B_n$ are present.
Conditioning on such an event is, by the construction of
$\mu_3$, the same as doing the electrical computations
on a contraction-deletion of $B_n'$.  Thus
$\mu_3 (B_n', \C)  ( \cdot \| D)$ 
is a mixture of $\mu_3 (G, \C') (\cdot )$ as $G$ ranges over 
contraction-deletions of $B_n'$ (where $\C'$ is what's left of
the enumeration when you get to $B_n$).  The claim is then
just Raleigh's monotonicity; $\mu_3 (G, \C') (C(A))$ is
a product of conditional probabilities $p_j$ as in the proof
of Theorem~\ref{convergence}; any contraction-deletion
of $B_n'$ can be contracted to $G_2$ or deleted to $G_1$;
monotonicity says that contracting increases total current
and deleting decreases it, so each $p_j$ increases with 
deletion and decreases with contraction, and the claim is shown.

It remains to show that $\mu_3 (G_1) (C(A)) - \mu_3 (G_2) (C(A))
\rightarrow 0$ as $n \rightarrow \infty$ for each $A$.  For this,
use the random walk scenario.  Let $B_M$ be a box containing $A$.
Let $\ee > 0$ be arbitrary and $L$ be large enough so that 
the union of $L$ independent SRW's started anywhere on $\partial B_M$
will cover all the edges of $A$ with probability at
least $1-\ee$.  The following fact can be found in or
deduced from Lawler (1991): the hitting measure
of $\partial B_M$ for SRW on $B_n$ started from the vertex
$v$ converges as $v$ goes to infinity and $n$ varies arbitrarily
with $v \in B_n$.  This implies that for sufficiently large $n$,
the total variation distance between the hitting measures
on $\partial B_M$ from any two vertices on $\partial B_n$ 
can be made less than $\ee / L$.  

Now view $G_1$ and $G_2$ as graphs and couple SRW's 
$\gamma_i$ from the origins on $G_i$ as follows.  They 
are the same until they hit the boundary (which has been collapsed
to a single point in $G_2$).  Then they are coupled so that 
their next hits of $\partial B_M$ occur in the same place
(though not necessarily at the same time) with probability
as close to one as possible; this probability is at least
$1 - \ee / L$.  Then they make the same 
moves until they hit $\partial G_i$, become recoupled as often as
possible when they hit $\partial B_M$ again, and so on.  
The probability is at least $1 - \ee$ that 
$\gamma_1$ and $\gamma_2$ are coupled whenever they
are inside $B_M$ up to the first $L$ hits of 
$\partial B_M$.  At this point, the probability is at least 
$1 - \ee$ that all edges in $B_M$ have been traversed,
in which case the subgraph $\Tree (\gamma_1 )$ is in
the event $C(A)$ if and only if $\Tree (\gamma_2 )$ is. 
Thus $|\mu_3 (G_1 ) (C(A)) - \mu_3 (G_2 )(C(A))| < 2\ee$. 
Since $\ee$ was arbitrary, that sandwiches $\mu (C(A) \|
C_n)$ between sequences with the same limit and proves the theorem.
$\Cox$.

\section{Loop-erased random walk} \label{LERW}

This section contains lemmas about loop-erased random walk.
The reason that loop-erased random walk is relevant to
this paper will be clear later but briefly it is
the following: when $\mu_2 (G,v)$ is used to
construct a random spanning tree on $G$, the unique
path connecting a vertex $w$ to $v$ is given by
a loop-erased random walk from $w$ to $v$.
The section is self-contained, but not formal.  For a more
complete development, see Lawler (1991; or 1980, 1983 and 1986).

Let $G$ be any graph and let $\gamma$ be a path on $G$.  The
following notational conventions will be used throughout.
The $i^{th}$ vertex visited by $\gamma$ is denoted $\gamma (i)$,
beginning at $\gamma (0)$.  If $\gamma$ is finite then 
$l (\gamma )$ denotes the length of $\gamma$ and
$\gamma'$ denotes the time reversal of $\gamma$, so 
$\gamma' (0) = \gamma (l(\gamma ))$.  If in
addition there is a path $\beta$ with $\beta (0) = \gamma' (0)$
then $\gamma * \beta$ denotes $\gamma$ followed by $\beta$.
The paths $\beta$ and $\gamma$ are said to intersect 
whenever $\beta (i) = \gamma (j)$ for some $i$ and $j$ not necessarily
equal but not both zero.  Finally, $\gamma \wedge n$ denotes 
the initial segment $\{ \gamma (i) : i \leq n \}$ of $\gamma$
and $\gamma \vee n$ denotes $\gamma$ from step $n$ onwards,
so $\gamma = (\gamma \wedge n) * (\gamma \vee n)$.

For finite paths $\gamma$ the loop-erasure operator $LE$ is
defined intuitively as follows.  If $\gamma$ is a self-avoiding
path (meaning that the vertices $\gamma (i)$ are distinct)
then $LE (\gamma ) = \gamma$.  Otherwise, the first time $\gamma$
visits a vertex $v$ twice, erase the loop at $v$.  In
other words, if $\gamma (i) = \gamma (j), \, i<j$ and $j$ is minimal
for this, delete from the sequence $\{ \gamma (k) \}$  
all the vertices with $i < j \leq k$.  If the result is still
not self-avoiding then repeat this step until it is.
The map $LE$ preserves the initial and final points of
a path.  For a given initial and final point $LE$ maps onto
the set of self-avoiding paths with the given endpoints but
is not one to one.  Let $\alpha$ be a self-avoiding path and
$m$ a positive integer and, following Lawler (1983) in slightly 
different notation, define $\Gamma^m (\alpha )$ to
be $LE^{-1} (\alpha ) \cap \{ \mbox{paths of length }m \}$.

If $\gamma$ is an infinite path that hits every vertex
finitely often then the paths $LE (\gamma \wedge n )$ converge
to an infinite path which will be called $LE (\gamma )$.
When $G = \Z^d$, $d \geq 3$ and $\gamma$ is a SRW from some
vertex $v$, then $\gamma$ hits each vertex finitely often
almost surely.  Consequently $LE (\gamma )$ is almost surely
well-defined.  The law of $LE (\gamma )$ is called 
the loop-erased random walk measure on $\Z^d$ from $v$, or
simply LERW.  (LERW can be defined on $Z^2$ as well but will
not be needed here.)

Commonly, an alternative construction for LERW is used.
Let $\gamma (0)$ be given and let the measure of the
event $\gamma (1) = v$ be given by the probability
the $\beta (1) = v$ where $\beta$ is a SRW conditioned
never to return to $\gamma (0)$.  In general, let the
measure of $\gamma (i+1) = v$ conditional on $\{ \gamma (j) :
j \leq i \}$ be given by the probability that 
$\beta (1) = v$ where $\beta$ is a SRW from $\gamma (i)$
conditioned never to return to $\{ \gamma (j), j \leq i \}$.
A similar construction gives the law of $LE (\gamma )$
when $\gamma$ is a SRW from $v$ on a finite graph $G$, stopped upon
hitting some vertex $w$.  In this case the conditional
probability of $\gamma (i+1) = v$ given $\gamma (j)$ for 
$j \leq i$ is given by the next step of a random walk 
conditioned to hit $w$ before returning to  $\{ \gamma (j), j \leq i
\}$.  These characterizations are easy to prove and
will be assumed freely when convenient.

\begin{lem} \label{dimension}
Let $v$ and $w$ be vertices in $\Z^d$, $d \geq 3$.  Let 
$\beta$ and $\gamma$ be independent LERW from $v$ and
SRW from $w$ respectively.  Then if $d = 3$ or $4$, $\beta$ and
$\gamma$ intersect infinitely often almost surely.
On the other hand if $d \geq 5$, $\beta$ and $\gamma$ 
intersect finitely often almost surely and the 
probability that they intersect at all (other than at $v$
if $v=w$) is bounded
between $c_1 (d) \, |v-w|^{4-d}$ and $c_2 (d) \, |v-w|^{4-d}$
for some constants $0 < c_1 (d)< c_2 (d)< \infty$.  
\end{lem}

\noindent{Proof:}  The statement for $d=3$ is proved in
Lawler (1988 equation 3.1) and for $d=4$ is proved in
Lawler (1986 Theorem 5.1).  For $d \geq 5$ the
fact that $\beta$ and $\gamma$ intersect finitely often
almost surely can be deduced from the corresponding facts
for two SRW's and the fact that LERW is a subsequence of
SRW.  To prove the quantitative bounds for $d \geq 5$, 
proceed as follows.

Let $X$ be the random number of intersection points of
a LERW from $v$ and an independent SRW from $w$, counted
with multiplicity $k$ if the point is hit $k$ times by the SRW.
The upper bound is a consequence of the following upper bound on 
$\E X^2$ which can be found in Lawler (1991 Chapter 3).
\begin{equation} \label{upper X^2}
\E X^2 \leq c |v-w|^{4-d} .
\end{equation}
Since $X$ is an integer-valued random variable, this immediately
establishes that $\P (X > 0) \leq c |v-w|^{4-d}$, which is
the desired upper bound on the probability that LERW from
$v$ intersects an independent SRW from $w$.  The lower bound
will be proved by showing
\begin{equation} \label{lower X}
\E X \geq c |v-w|^{4-d} .
\end{equation}
To see that~(\ref{upper X^2}) and~(\ref{lower X}) actually
imply $\P (X > 0) \geq c |v-w|^{4-d}$, write
\begin{eqnarray*}
\E X^2 & = & \P (X > 0) \E (X^2 \| X > 0) \\[2ex]
& \geq & \P (X > 0) (\E (X \| X > 0))^2 \\[2ex]
& = & \P (X > 0) \left [ {\E X \over \P (X > 0)} \right ]^2 \\[2ex]
& = & (\E X)^2 \P (X > 0 )^{-1} ,
\end{eqnarray*}
hence $\P (X > 0) \geq (\E X)^2 / \E X^2 $.

To show (\ref{lower X}), let $\beta$ be a SRW from $v$ and 
$\gamma = LE (\beta)$ be the corresponding LERW from $v$.
Write $G(x,y)$ for the Green's function,
i.e. the expected number of visits to $y$ of SRW starting at $x$.
It is known (e.g. Lawler 1991) that $G (x,y)$ is bounded between 
constant multiples of $|x-y|^{2-d}$ in each dimension $\geq 3$;
in this regard, let $0^{-n}$ denote the constant $G(x,x)$ to
avoid making explicit exceptions for zero in the summations.
Then~(\ref{lower X}) is implied by
\begin{equation} \label{lower LERW}
\P (x \in \gamma ) \geq c |v-x|^{2-d}
\end{equation}
since this implies
\begin{eqnarray*}
\E X & = & \sum_x \P (x \in \gamma ) G(w,x) \\[2ex]
& \geq & c\sum_s |\{ x: |v-x| = s \}| s^{2-d} (s+|v-w|)^{2-d} \\[2ex]
& \geq & c \sum_s s^{2-d} s^{d-1} (s+|v-w|)^{2-d} \\[2ex]
& \geq & c \sum_{s \geq |v-w|} s^{2-d} s^{d-1} (2s)^{2-d}
\end{eqnarray*}
which is just $c|v-w|^{4-d}$.

Finally, to show~(\ref{lower LERW}) let $\tau$ be the first
time (possibly infinity) that $\beta$ hits $x$ and write
$$\P (x \in \gamma ) \geq \P (\tau < \infty ) \P (\beta \wedge
\tau \mbox{ is disjoint from } \beta \vee \tau   
\| \tau < \infty) ,$$
The first factor is at least $c|v-x|^{2-d}$ so it remains
to bound the second factor away from zero.  Since $\beta \vee \tau$ 
is independent of $\beta \wedge \tau$ given $\tau < \infty$,
the second factor is the probability that $\beta \wedge \tau$
is disjoint from an independent SRW $\beta_1$ from $x$,
where $\beta$ is a SRW from $v$ conditioned to hit $x$.
Write $\beta_2 = (\beta \wedge \tau )'$, 
so $\beta_2$ is a SWR from x conditioned to hit $v$ 
and stopped when it hits $v$.  Since two independent SRW's from
$x$ are disjoint with positive probability for $d \geq 5$, 
it remains to show that conditioning one of the walks to reach $v$
does not alter this.  We may assume that $|v-x|$ is greater than
some fixed constant $r_0$, since~(\ref{lower LERW}) is immediate
for $|v-x| \leq r_0$ just from transience of the SRW.  

Let $\gamma_1$ and $\gamma_2$ be independent SRW's from $x$.
Fix any positive $\ee$.
Since independent SRW's from $x$ intersect finitely often with
probability one, an $M$ can be chosen large enough so that 
$\P(\gamma_1 \vee M \mbox{ intersects } \gamma_2) < \ee$.
By transience of SRW, an $M' > M$ can be chosen so that
$\P(\gamma_2 \vee M' \mbox{ intersects } B(x,M)) < \ee$,
where $B(y,k)$ is the cube of radius $k$ centered at $y$.
It is known, via triviality of the Martin boundary for SRW
(e.g. Lawler 1991 Chapter 2), that SRW from
$x$ conditioned to hit $y$ converges weakly to unconditioned
SRW from $x$ as $|x-y| \rightarrow \infty$, so $r_0$ may
be chosen such that $|x-y| \geq r_0/4$ implies that the
total variation difference between $\gamma_1 \wedge M$
and $\beta_2 \wedge M$ is less than $\ee$.  Similarly,
let $r = |x-v|$ and let $\alpha$ be a SRW from $x$ conditioned
to avoid $B(v, 3r/4)$; then the same argument about the
Martin boundary shows that the distribution of $\alpha$ 
converges weakly to that of $\gamma_2$ as $r \rightarrow \infty$,
so $r_0$ can be chosen large enough so that $r \geq r_0$
implies that the total variation distance between $\gamma_2
\wedge M'$ and $\alpha \wedge M'$ is at most $\ee$.

Now let $p_1 = \P (\gamma_1 \wedge M \mbox{ is disjoint from }
\gamma_2)$.  Let $p_2 = \P (\gamma_2 \mbox{ is disjoint from }
B(v, 3r/4))$ and let $p_3 = \min_{y \in \partial B(v,r/2)}
\P$ (SRW from $y$ conditioned to hit $v$ does so before
leaving $B(v,3r/4))$.  Note that $p_1$ is bounded away from zero
by the standard result, while $p_2$ and $p_3$ are easily seen 
by scaling to be bounded away from zero in any fixed dimension.
Let $\sigma$ be the first time $\beta_2$
hits $B(v,r/2)$ and write 
\begin{eqnarray*}
&& \P (\beta_2 \dj \beta_1) \\[2ex]
& \geq & p_2 \P (\beta_2 \dj \alpha ) \\[2ex]
& \geq & p_2 \P (\beta_2 \wedge \sigma \dj \alpha) \P( \beta_2
   \vee \sigma \dj \alpha \| \alpha ) \\[2ex]
& \geq & p_2 [ \P (\beta_2 \wedge M \dj \alpha \wedge M') - 
   \P (\beta_2 \wedge M \mbox{ intersects } \alpha \vee M')\\
& & - \P ((\beta_2 \wedge \sigma ) \vee M \mbox{ intersects } 
   \alpha )] \P (\beta_2 \vee \sigma \dj \alpha \| \alpha ) \\[2ex]
& \geq & p_2 [\P (\gamma_1 \wedge M \dj \gamma_2 \wedge M') - 2\ee
   - \P (\alpha \vee M' \dj B(x,M)) \\
&&   - \P ((\beta_2 \wedge \sigma) 
   \vee M \mbox{ intersects } \alpha)] p_3 \\[2ex]
& \geq & p_2 p_3 [ p_1 - 2\ee - 2\ee - \P ((\beta_2 \wedge \sigma)
   \vee M \mbox{ intersects } \alpha)] 
\end{eqnarray*}
by choice of $M$ and $M'$.  Since $p_i$ are all bounded away from
zero, it remains to show that $\P ((\beta_2 \wedge \sigma)
\vee M \mbox{ intersects } \alpha)$ is small.  But the
distribution of $\beta_2 \wedge \sigma$ is given by
a SRW conditioned to hit $B(v,r/2)$ at some random
point $y$, stopped when it does so, reweighted by 
$\P(\mbox{SRW from $y$ hits } v)$ and normalized.  Scaling shows
that $\P (\mbox{SRW from $x$ hits } B(v,r/2))$ is
bounded below, and as $y$ varies over the boundary of $B(v,r/2)$ 
in a fixed dimension, the ratios of these reweights are bounded.  
Thus the Radon-Nikodym derivative
${d(\beta_2 \wedge \sigma) \over d( SRW \wedge \sigma)}$ 
is bounded above, and hence $\P ((\beta_2 \wedge \sigma) \vee M
\mbox{ intersects } \alpha$ is bounded by a constant times
$\P (\gamma_1 \vee M \mbox{ intersects } \alpha$ and
the latter is at most $p_2^{-1} \ee$.  This completes the
proof that $\P (\beta_2 \dj \beta_1)$ is bounded away from
zero, thus proving~(\ref{lower LERW}) and~(\ref{lower X}).   $\Cox$.

\begin{lem} \label{reversal}
Let $G$ be any graph and $\alpha$ a finite path in $G$.  
Let $\Phi^m (\alpha ) = \{ \beta : \beta' \in \Gamma^m 
(\alpha' ) \}$ be the set paths of length $m$ whose
``backward loop-erasure'' is $\alpha$.  
Then for each $m$ there is a bijection $T^{m,\alpha}$
between $\Gamma^m (\alpha )$ and $\Phi^m (\alpha )$
such that the multiset of sites visited by $\gamma$
is the same as the multiset of sites visited by $T^{m,\alpha}
(\gamma )$.
\end{lem}

\noindent{Proof:}  Lawler (1983 Proposition 2.1) states this
for $G = \Z^d$ and for sets instead of multisets.  The 
proof actually shows that multisets are preserved.  Clearly,
if the proposition is true for $\Z^d$ it is true for
subgraphs of $\Z^d$, which is all that is used below.
It is easy, however, to see that Lawler's proof is valid for
any graph.     $\Cox$

\begin{lem} \label{firstMsteps}
Let $w$ be any vertex in $\Z^d$, $d \geq 3$.  For any
positive integer $L$, let $x$ be a vertex in $B_n$ at distance
at least $L$ from $w$, where $B_n$ is large enough to
contain $w$.  Let $\gamma$ be a SRW from $x$ on $B_n$ 
conditioned to hit $w$ before returning to $x$ and let 
$\alpha = LE (\gamma')$.  Then 
the distribution of the first $M$ steps of $\alpha$ 
converges as $n,L \rightarrow \infty$ to the distribution 
of the first $M$ steps of LERW on $\Z^d$ from $w$,
the convergence being uniform over choices of $x$.
\end{lem}

\noindent{Proof:}  First note that by time reversal, $\gamma'$ 
is distributed as SRW from $w$ conditioned to hit $x$ before
returning to $w$.   It suffices to show that for each 
self avoiding path $\beta$ of length $j < M$ from $w$, and
each neighbor $v$ of $\beta (j)$,
the conditional probability that $\alpha \wedge {j+1} (j+1) = v$ 
given $\alpha \wedge j = \beta$ approaches 
$\P (LERW \wedge {j+1} (j+1) = v \| LERW \wedge j = \beta )$.
By the alternative construction for LERW, the latter probabilities
for fixed $\beta$
are proportional to the quantities $p (v)$ defined by $p(v) = 
\P (SRW \mbox{ from } 
v \mbox{ never hits } \beta )$,
and are thus given by $p (v)$ normalized to sum to one.
Similarly, the former probabilities are proportional
to $q(v) = \P (\mbox{SRW on $B_n$ from } 
v \mbox{ hits $x$ before hitting } \beta)$.  

Let $K$ be the box such that $x \in \partial B_K$, so for
fixed $w$, $K \rightarrow \infty$ as $L \rightarrow \infty$.
Let $Q( \cdot )$ be the hitting measure on the boundary of
$B_K$ for SRW from $w$ conditioned to avoid $\beta$. 
It is known (e.g. Lawler 1991 Theorem 2.1.2) that $Q(y)$
is bounded between $1 - \ee (K, \beta )$ and $1 + \ee (K, \beta 
)$ times the hitting measure for SRW starting from the origin,
where $\ee (K, \beta ) \rightarrow 0$ as $K \rightarrow \infty$.
To make use of this, write
\begin{eqnarray}
q(v) & = & \P (\mbox{SRW on $B_n$ from $v$ hits the boundary of }
    B_K \mbox{ before hitting } \beta \label{avg} \\[2ex]
& \times & \sum_{y \in \partial B_K} Q(y) \P (\mbox{SRW on }
    B_n \mbox{ from $y$ hits $x$ before } \mbox{hitting}
    \beta (0) , \ldots , \beta (j)) . \nonumber 
\end{eqnarray}
The first factor on the RHS of equation~(\ref{avg}) converges
to $p(v)$ as $K \rightarrow \infty$.  The second
one, according to the observation about $Q$ above, 
may only vary with $v$ by a factor of at most $1 \pm
\ee (K, \beta )$.  Thus for fixed $\beta$, 
$q(v)$ normalized converges to $p(v)$ normalized as
$n,K \rightarrow \infty$ uniformly in $x \in \partial B_K$,
hence as $n,L \rightarrow \infty$ uniformly in $x$ at
distance at least $L$ from $y$, and the proof is done.    $\Cox$

\begin{lem} \label{strong reversal}
Remove the conditioning in Lemma~\ref{firstMsteps} so that
$\gamma$ may return any number of times to $x$ before hitting
$w$.  Then (i) the conclusion that $\alpha \wedge M$ converges
to $LERW |M$ uniformly in $x$ still holds; (ii) $(LE (\gamma'))'$ 
has the same distribution as $LE (\gamma )$.
\end{lem}

\noindent{Proof:}  For finite paths $\beta$ from $x$ in $B_n$,
let $W(\beta ) = W(B_n , \beta )$ denote $\P (\gamma \wedge {l( \beta )}
= \beta )$, which can be written as $\prod_i (\mbox{number of
neighbors of } \beta (i))^{-1}$.  To prove $(ii)$, write
\begin{equation} \label{path measures}
\P (LE (\gamma' ) = \alpha') = \sum_m \sum_{\beta \in \Gamma^m 
    (\alpha' ) \cap S} W(\beta' )
\end{equation}
where $S$ is the set of paths that never return to $w$.
Since the bijections $T^{m,\alpha}$ of Lemma~\ref{reversal}
preserve the multiset of sites visited, they preserve $W$
and can be used to rewrite~(\ref{path measures}) as
$$ \sum_m \sum_{\beta \in \Phi^m (\alpha' ) \cap S} W(\beta' ) $$
which is by definition of $\Phi^m$ just
$$ \sum_m \sum_{\beta' \in \Gamma^m (\alpha ) \cap S} W(\beta' ) $$
which is $\P (LE (\gamma ) = \alpha )$.

To prove (i) note that the distribution of $LE (\gamma )$
is independent of the number of times $\gamma$ returns to $x$.
Then by $(ii)$, the distribution of $LE (\gamma' )$ is 
independent of the number of times $\gamma$ returns to $x$.
In particular it is unaffected by conditioning on this number
being zero, thus Lemma~\ref{firstMsteps} holds even after
conditioning.   $\Cox$

\section{Number and shape of the components} \label{shape}

The following easy lemma connects loop-erased random walk
to the random walk method of generating a random spanning
tree of a finite graph.  Recall the definition of $\Tree (\gamma )$
at the beginning of section 2.  

\begin{lem} \label{path=LE}
Let $v$ and $w$ be distinct vertices of a finite graph $G$
and let $\gamma$ be any path from $v$ to $w$, not necessarily
self-avoiding.  Then the unique path connecting $w$ to
$v$ in $\Tree (\gamma )$ is given by $LE (\gamma' )$.
\end{lem}

\noindent{Proof:}  Let $\alpha$ be $LE (\gamma' )$ and $\beta$
be the path connecting $w$ to $v$ in $\Tree (\gamma )$.
Clearly $\alpha (0) = \beta (0) = w$.  Now assume for induction
that $\alpha (i) = \beta (i)$.  Then $\beta (i+1)$
is the unique $x$ for which $\Tree (\gamma )$ has an
oriented edge $\overline {x \, \beta (i)}$.  This is just
$\gamma (j-1)$ where $j$ is minimal such that $\gamma (j) = 
\beta (i)$.  This is also equal to $\gamma' (j+1)$ where $j$ 
is maximal for $\gamma' (j) = \beta (i) = \alpha (i)$.  Then
when applying loop-erasure to $\gamma'$, the edge from 
$\alpha (i)$ to $x$ is never erased, hence $\alpha (i+1) = x$.
By induction, $\alpha = \beta$.   $\Cox$

The main theorem on connectedness can now be proved.
\begin{th} \label{tree-forest}
Let $\mu$ be the limiting measure on subgraphs of $\Z^d$,
$d \geq 3$ constructed in section~\ref{equivalences}.
Then for $d = 3$ or $4$, $\mu$ concentrates on connected
graphs.  For $d \geq 5$, $\mu$ concentrates on graphs
with infinitely many components.  In this case, 
$|v-w|^{d-4} \P (v \mbox{ and } w \mbox{ are } \\ \mbox{connected})$
is bounded between $c_1 (d)$ and $c_2 (d)$ for 
$v \neq w$ and some constants $0 < c_1 (d) < c_2 (d) < \infty$.
\end{th}

\noindent{Proof:}  Fix $d$ for the moment.  If $d=2$,
$\mu$ can be defined via $\mu_2$ without a limiting 
procedure, since SRW in $Z^2$ hits every point, and connectedness
follows immediately.  So assume
without loss of generality that $d > 2$.   Let $v$ and $w$ be
distinct vertices.  The main project will be determining whether
$v$ and $w$ are almost surely connected.  If so, then by
countable additivity the whole graph is almost surely connected.
If not, then another few sentences will show that there
are almost surely infinitely many components.

Fix the vertices $v$ and $w$.
The argument will use the random walk scenario, writing
$\mu$ as the limit of $\mu^n = \mu_2 (B_n , v)$ as $n \rightarrow
\infty$.  Let $C$ be the event $\{ v \mbox{ is connected to }w \}$.
Since the convergence is weak, and the indicator function
$\one_C$ is not continuous, $\mu^n (C)$, which is
always $1$, does not necessarily converge to $\mu (C)$.  	
To get information about $\mu$ we must work instead with the
continuous events $C_M = \{ v \mbox{ is connected to } w
\mbox{ by a path of length } \leq M \}$.  Specifically,
weak convergence implies $\mu^n \rightarrow \mu$
on each $C_M$, hence 
\begin{equation} \label{limit form}
\mu (C) = \lim_{M \rightarrow  \infty} \lim_{n \rightarrow 
    \infty} \mu^n (C_M) .
\end{equation}
Another way to say this is to let $L_n$ be the length of
the path connecting $v$ and $w$ under $\mu^n$.  Then 
$v$ and $w$ are $\mu$-almost surely connected if and
only if the $L_n$'s are tight.  Equation~(\ref{limit form}) will be
used to show that $\mu (C)$ is equal to the probability that LERW
from $w$ intersects an independent SRW from $v$ 
(equation~\ref{goal} below).

To analyze $\mu^n$, run a SRW $\beta$ from $v$ on $B_n$.
Let $\tau$ be the first time $\beta$ hits $w$ and let
$\gamma = \beta \wedge \tau$.  The path connecting $v$
and $w$ in $\Tree (\beta )$ is determined by $\gamma$.  There
are two possibilities: either $\beta$ hits $\partial B_n$ before
hitting $w$ or vice versa.  If it hits $w$ first, it is easy
to check that the conditional distribution of the 
length of $\gamma$ is tight as $n \rightarrow \infty$.  

To examine the other possibility, condition (hereafter)
on $\beta$ hitting $\partial B_n$ before $w$ and let $x$ be
the first point where $\beta$ hits $\partial B_n$.
Write $\gamma = \gamma_1 * \gamma_2$ where $\gamma_1$
is the initial segment of $\gamma$ up to the first hit
of $x$ and $\gamma_2$ is all the rest.  Then $\gamma_1$
is distributed as SRW from $v$ stopped upon hitting
the boundary and conditioned to do this before it hits
$w$.  Then as $n \rightarrow \infty$ the first $M$ steps
of $\gamma_1$ converge for each $M$ to the first $M$ steps
of an infinite SRW from $v$ conditioned never to hit $w$.

Recall from Lemma~\ref{path=LE} that the path connecting 
$w$ to $v$ is given by $LE (\gamma' ) = 
LE (\gamma_2' * \gamma_1' )$.  Fix any $M$.
Observe that $LE (\gamma_2' ) \wedge M = LE (\gamma' ) \wedge M$
whenever $LE (\gamma_2' ) \wedge M$ is disjoint from $\gamma_1$.
This is because $LE (\gamma' ) = LE (LE (\gamma_2' ) *
\gamma_1' )$ and the addition of $\gamma_1'$ cannot alter
any initial segment of $LE (\gamma_2' )$ that it does not intersect.
It should now be clear where Lemma~\ref{dimension} comes in;
the rest of the work will be in identifying the distributions
of $\gamma_1$ and $LE (\gamma_2' )$ and taking limits
correctly.

Let $\alpha$ be a LERW from $w$ independent from $\beta$.
Recall from Lemma~\ref{strong reversal} that 
\begin{equation} \label{pair conv}
LE (\gamma_2 ') \wedge M \dconv \alpha \wedge M 
\end{equation}
as $n \rightarrow \infty$,
even when conditioned on $x$.  (Here the dependence of $\gamma_2$ on
$n$ is supressed in the notation.)  Since $\gamma_1$
and $\gamma_2$ are conditionally independent given $x$,
it follows that for any $M$, the pair $(LE (\gamma_2' ) 
\wedge M, \gamma_1 \wedge M)$ converges to
$(\alpha \wedge M, \beta \wedge M)$ as $n \rightarrow \infty$.

Let $D$ be the event that $\alpha$ and $\beta$ intersect.  
Let $D_M$ be the event that $\alpha \wedge M$ and $\beta \wedge M$ 
intersect, and let $D_M'$ be the event that $\alpha \wedge M$ 
and $\beta$ intersect.  Then $D_M , D_M' \uparrow D$, so
$\P (D_M ) , \P (D_M') \uparrow \P (D)$.

Recall that $C_M$ is the event that the path connecting 
$v$ to $w$ in $\Tree$ has length at most $M$.  Then 
\begin{equation} \label{sandwich}
\gamma_2' \wedge M \cap \gamma_1 \wedge M \neq \emptyset 
    \Rightarrow C_{2M} \Rightarrow \gamma_2' \wedge 2M 
    \cap \gamma_1 \neq \emptyset .
\end{equation}
It follows from~(\ref{pair conv}) that
$$\lim_{n \rightarrow \infty} \mu^n (LE (\gamma_2') \wedge M 
    \cap \gamma_1 \wedge M \neq \emptyset ) = \P (D_M) .$$
Let $u (M)$ be large enough so that $\P (\alpha \wedge 2M
\cap \beta \wedge u(M) \neq \emptyset \| \alpha \wedge 2M
\cap \beta \neq \emptyset ) > 1 - 1/M$.
Then it also follows from~(\ref{pair conv}) that
$$\lim_{n \rightarrow \infty} \mu^n (LE (\gamma_2') \wedge 2M 
    \cap \gamma_1 \neq \emptyset ) \leq {M \over M-1} \P 
    (\alpha \wedge 2M \cap \beta \wedge u(M) \neq
    \emptyset \leq {M \over M-1} \P (D_{2M}') .$$
Now taking limits as $n \rightarrow \infty$ of~(\ref{sandwich})
gives 
$$ \P (D_M) \leq \lim_{n \rightarrow \infty} \mu^n (C_{2M}) \leq
    {M \over M-1} \P (D_{2M}') .$$
Taking the limit in $M$ and using equation~(\ref{limit form}) gives
\begin{equation} \label{goal}
\P (D) = \mu (C) . 
\end{equation}
Now if $d = 3$ or $4$, Lemma~\ref{dimension} says that 
the probability of $\alpha$ intersecting an independent SRW from
$v$ is one; since $\beta$ is distributed as an independent SRW
from $v$ conditioned on an event of positive probability, this means
$\P (D) = 1$, from which the statement of the theorem follows
immediately.  

On the other hand, consider the case $d \geq 5$.  By 
Lemma~\ref{dimension},
the probability that $\alpha$ intersects an independent $SRW$ 
from $v$ is bounded between constants times $|v-w|^{4-d}$.
Since the event that SRW from $v$ actually hits $w$
is of order $|v-w|^{2-d}$, $\beta$ is distributed as a SRW
conditioned on an event of probability $1-c|v-w|^{2-d}$, and
it follows from $\P(A) / \P(B) \geq \P ( A \| B) \geq
(\P (A) - \P(B^c)) / \P (B)$ that $\P (D)$ 
is bounded between constant multiples of $|v-w|^{4-d}$,
hence $\P (C)$ is also, which was to be shown.
It follows immediately that the measure $\mu$ does 
not concentrate on connected graphs.  

To see that the measure concentrates on graphs with infinitely
many components, recall from Theorem~\ref{convergence} that
$\mu$ is stationary and from Theorem~\ref{tail} that the
tail field is trivial.  Then $\mu$ is ergodic, so the number
of components is some constant $K$ almost surely.  
To bound $K$, write $I(x,y)$ for the indicator
function of the event that $x$ is connected to $y$ and calculate 
$$\E \sum_{x,y \in B_n} I(x,y) = \sum_{x,y \in B_n} 
\E I(x,y) = \sum_{x,y \in B_n} O(|x-y|^{4-d}) = O(n^{d+4}) .$$
On the other hand, if $B_n$ is partitioned into at most 
$K$ connected components, $K < \infty$, then
$$\sum_{x,y \in B_n} I(x,y) \geq n^{2d}/K .$$
When $d \geq 5$ this is greater than $O(n^{d+4})$ for any 
finite $K$, so $K$ must be infinite almost surely.
This completes the proof.    $\Cox$

The last theorem is about the shape of the tree when
$(d \leq 4)$.

\begin{th} \label{one-end tree}
If $d \leq 4$ then the measure $\mu$ concentrates on trees
having only one topological end, i.e. trees for which
removal of any vertex divides the tree into components
precisely one of which is infinite.
\end{th}

\noindent{Proof:}  Call a vertex $x$ in a subgraph of
$\Z^d$ a separator if removal of $x$ leaves more than 
one infinite component.  Call $x$ a branchpoint if its
removal leaves more than two infinite components.  
Burton and Keane (1989 Theorem 2) show that the set
of branchpoints for a subgraph of the integer lattice may not be a set
of vertices of positive density.  By stationarity and
ergodicity it follows that there are no branchpoints
at all almost surely.  Then the tree has at most two 
topological ends.  

The number of topological ends is translation invariant,
hence almost surely constant.  Assume for contradiction 
that there are almost surely two.  Then the spanning
tree $\Tree$ looks like a doubly infinite line to which
has been attached at each vertex a finite (possibly empty) 
tree.  The vertices on the infinite line are precisely
those vertices that are separators and by ergodicity
and tail triviality this set has a density $D_{sep} > 0$ that
is almost surely constant.  

For any vertices $v_1, v_2$ and $v_3$ say that $v_2$
separates $v_1$ and $v_3$ if the unique path in $\Tree$
from $v_1$ to $v_3$ passes through $v_2$.  Observe
that if $v_1, v_2$ and $v_3$ are all on the infinite
line in $\Tree$ then one of them separates the other two.
Thus for any $v_1 , v_2, v_3$,
$$\sum_i \P (v_i \mbox{ separates the other two})
    \geq \P (v_1 , v_2 , v_3 \mbox{ are all separators}) .$$
Now triviality of the tail implies that $\mu$ is mixing
of all orders, and in particular $2$-mixing implies
$$\P (v_1 , v_2 , v_3 \mbox{ are all separators}) 
    \rightarrow D_{sep}^3 $$
as the pairwise distances $|v_i - v_j|$ all go to infinity.
To get a contradiction then, it suffices to show that
\begin{equation} \label{separation}
\P (x \mbox{ separates } v \mbox { and } w)
   \rightarrow 0
\end{equation}
as the pairwise distances between $v,w$ and $x$ all go to infinity.

Assume then that the pairwise distances between the vertices
$v,w$ and $x$ are at least $L$ for some $L > 0$.  Use
the random walk scenario with $\mu = \lim \mu^n$ where
$\mu^n$ is constructed as $\mu_2 (B_n , v)$ for $B_n$ large
enough to contain $v,w$ and $x$.  Fix $n$ for the moment
and let $\gamma$ be the initial
segment of the random walk from $v$ up to the first hitting
of $w$.  Here is how $\gamma$ determines whether $x$ separates
$v$ and $w$.  If $\gamma$ does not hit $x$ then $x$ does
not separate $v$ from $w$.  If $\gamma$ does hit $x$ then
let $\gamma_1$ be $\gamma$ up to the first hitting of $x$
and $\gamma_2$ be the rest of $\gamma$.  The path connecting
$v$ and $w$ in $\Tree$ is the path connecting them in
$\Tree (\gamma )$ which is given by 
$$LE (\gamma') = LE (\gamma_2' * \gamma_1') 
    = LE (LE (\gamma_2' ) * \gamma_1' ) .$$
Now $x$ appears only once in $LE (\gamma_2') * \gamma_1'$,
namely at the point where they join.  Thus $x$ separates
$v$ and $w$ if and only if it does not get erased when 
$LE$ is applied to $LE (\gamma_2') * \gamma_1'$.  If  
$\gamma_1'$ is disjoint from $LE (\gamma_2 )$ except at $x$,
it is clear that the loop-erasure on
$LE (\gamma_2') * \gamma_1'$ acts only on the $\gamma_1'$ part
and $x$ never gets erased.  Conversely, the first time that
$\gamma_1'$ intersects $LE (\gamma_2')$, the vertex $x$
will be erased.  Therefore, $x$ is erased if and only if 
$\gamma_1'$ and $LE (\gamma_2')$ are disjoint except
at $x$.  It remains to show that the probability of
these paths being disjoint goes to zero as $n \rightarrow \infty$
and then $L \rightarrow \infty$.   

For each $M$, the probability that $LE (\gamma_2')'\wedge M$ and 
$\gamma_1' \wedge M$ are disjoint are an upper bound for the
probability that $LE (\gamma_2')'$ and $\gamma_1'$ are 
disjoint.  Then to show~(\ref{separation}) it suffices to show:
\begin{equation} \label{separation2}
\inf_M \lim_{L \rightarrow \infty} \lim_{n \rightarrow \infty} 
    \P (LE (\gamma_2' )'\wedge M \cap \gamma_1' \wedge M \neq \{ x \})
    = 0 .
\end{equation}
Now by Lemma~\ref{strong reversal}~$(ii)$, $LE (\gamma_2')'$
has the same distribution as $LE (\gamma_2 )$.  Combine the 
fact that $\gamma_1$ and $\gamma_2$ are independent
with the fact from Lemma~\ref{firstMsteps} that $LE (\gamma_1)\wedge M$
converges to LERW$\wedge M$ and the fact that $\gamma_2 \wedge M$
converges to an independent SRW$\wedge M$ to 
rewrite~(\ref{separation2}) as
$$\inf_M \P (LERW\wedge M  \cap SRW\wedge M  \neq \{ x \} ) = 0,$$
where LERW and SRW are independent starting from $x$.  
This is a direct consequence of Lemma~\ref{dimension}.
Thus~(\ref{separation2}) and~(\ref{separation}) are shown
and the theorem is proved.   $\Cox$

\renewcommand{\baselinestretch}{1.0}\large

\end{document}